\documentclass[11pt]{article}
\usepackage{amssymb}
\usepackage[dvips, usenames]{color}
\usepackage{amsmath,amsfonts,amsthm,amscd,latexsym}
\def\uuar{\mathord{\mbox{\makebox[0pt][l]{\raisebox{.4ex}
{$\uparrow$}}$\uparrow$}}}

\def\dda{\mathord{\mbox{\makebox[0pt][l]{\raisebox{-.4ex}
{$\downarrow$}}$\downarrow$}}}

\def\da{\mathord{\downarrow}}

\def\ur{(U,R)}

\newtheorem{tm}{Theorem}[section]
\newtheorem{pn}[tm]{Proposition}
\newtheorem{lm}[tm]{Lemma}
\newtheorem{cy}[tm]{Corollary}

\theoremstyle{definition}
\newtheorem{dn}[tm]{Definition}
\newtheorem{rk}[tm]{Remark}
\newtheorem{ex}[tm]{Example}

\usepackage[a4paper]{geometry}
\allowdisplaybreaks[4]
\begin{document}
\title{When do CF-approximation spaces capture sL-domains
\thanks{Supported by the Natural Science Foundation of China (12231007, 12371462, 11671008), the Natural Science Foundation of Jiangsu Province (BK20170483)}}
\author{{Guojun Wu$^{1}$\thanks{Corresponding author.},\  \ Luoshan Xu$^{2}$,\ \   Wei Yao$^{1}$}\\
{\small  }\\
{\small $^{1}$School of Mathematics and Statistics}\\
{\small  Nanjing University of Information Science and Technology}\\
{\small Nanjing, 210044,  China}\\
{\small $^{2}$College of Mathematical Science, Yangzhou
University}\\
 {\small Yangzhou  225002,  China}}

 \date{}
\maketitle

\begin{abstract}
In this paper, by means of   upper approximation operators in  rough set theory, we study representations for sL-domains and its special subclasses.  We introduce the concepts of  sL-approximation
spaces, L-approximation
spaces and bc-approximation spaces, which are special types of CF-approximation
spaces. We prove that the collection of CF-closed sets in an sL-approximation
space (resp., an L-approximation
space, a bc-approximation space) ordered  by set-theoretic inclusion is an sL-domain (resp., an L-domain, a bc-domain);
 conversely, every sL-domain (resp., L-domain, bc-domain) is order-isomorphic to the  collection of CF-closed sets of an sL-approximation
space (resp., an L-approximation
space, a bc-approximation space). Consequently, we establish  an equivalence between  the category of sL-domains (resp.,   L-domains) with Scott continuous mappings  and  that of  sL-approximation spaces (resp.,  L-approximation spaces) with CF-approximable relations.\\

\noindent {\bf Keywords}\indent upper approximation operator;\ sL-domain;\ sL-approximation space;\\
\indent \indent \indent \indent\
\ L-domain;\ bc-domain; \ categorical equivalence
\end{abstract}

\section{Introduction}
Domain theory \cite{Gierz} is a branch of mathematics that deals with a specific class of partially ordered sets.
The domain structure serves as a denotational model for functional languages, playing a crucial role in theoretical computer science.
The research of domain theory can be   traced back to the pioneering work of  Scott \cite{Scott-72} in the  1970s. The basic feature of domain theory is the interaction between ordered structures and topological structures \cite{Jean-2013}. Recent academic focus has increasingly shifted towards exploring the nexus between domain theory and the mathematical foundations of artificial intelligence. This includes areas such as rough set theory \cite{Lei-Luo,Pawlak,Zou-Li-Ho}, formal concept analysis \cite{Ganter,Wsw-Lqg,Zxn-Wlc-Lqg}, non-classical logic \cite{Wang-log2,Wang-log1,Wang-log-2024}, and information systems \cite{He-Xu,Spreen-xu-inf,Wu-Guo-Li2,Wang-Li}.
The connections between domain theory and various mathematical structures are mainly reflected  in the representation/characterization of various domains \cite{wang-Zhou-Li-uni,XU-mao-form}. The representation of domains refers to that by utilizing an appropriate family induced by a specific mathematical structure, one can effectively  capture and reconstruct various domain structures, where this family is ordered by the set-theoretic inclusion \cite{Zou-Wang}.

Rough set theory \cite{Pawlak}, introduced by Pawlak and further developed by numerous mathematicians, holds fundamental significance  in artificial intelligence  in and computer science.
The core of rough set theory are  approximate spaces,  which are key mathematical structures.
A classical approximation space is essentially a  pair composed of a finite background set and an equivalence relation.   Subsequently,  by replacing equivalence relations with   binary relations, the notion of   approximation spaces was generalized to the so-called  generalized approximation spaces,  permitting non-finite background sets \cite{Kondo,Khan,gl-liu}.
The two approximation operators  defined by the binary relation of a generalized approximation space are similar to the topological inner operator and toplogical closure operator, which  are  mutually dual \cite{Yang-Xu-Top}.
J\"{a}rvinen, in \cite{jar2}, provided
the lattice-theoretical background and ordered properties  of generalized approximation space. In \cite{Lei-Luo}, Lei and Luo explored  representations of complete lattices and algebraic lattices based on approximation operators. Yang and Xu, in \cite{Yang-Xu-Alg}, delved into the ordered structures of various families of subsets induced by approximation operators. Zou and Li, in  \cite{Zou-Li-Ho}, characterized continuous posets and Scott closed sets by help with  approximation operators on posets. P. Eklund, et al.,  in  \cite{Eklund}, showed that partially ordered monads contain  appropriate structures for modeling generalized approximation spaces.  These studies reveal that there
are close links between generalized approximation spaces/rough set theroy and ordered structures/domain theory.

 Abstract bases, known for their capacity to generate various domains using round ideals \cite{xu-mao-ab}, possess a narrow research framework, potentially hindering connections with domain structures and other mathematical structures.
 To address this limitation, we want to extend the class of abstract base   to a broader class, ensuring that every object in the resulting larger class can still yield a continuous domain.
This extension adheres to the principle   that  the resulting  class should not be too broad, as it may lose connections with structures rooted in specific mathematical backgrounds, including topology,  mathematical logic, rough set theory, etc.
Notably, an abstract base $(B, \prec)$ is a  type of  generalized approximation space  from the perspective  of rough set theory. It is thus natural to wonder  whether one can appropriately extend abstract bases to a larger class of approximation spaces.
 Recently,
Wu and Xu in \cite{Wu-Xu} extended  an abstract basis to a CF-approximation
space  and generalized round ideals  to CF-closed sets. This extension demonstrated  that the class of abstract basis is a proper subclass of    CF-approximation
spaces.  The collection of CF-closed sets of a CF-approximation
spaces can generate a  domain and every  domain can   be reconstructed  from  a CF-approximation space (may not be a abstract base). This work gave a positive answer to the above problem to some extent.

The concept  of L-domains was introduced by  Jung  \cite{Jung1,Jung3}. It was proved that the category of L-domains (resp., algebraic L-domains) is a maximal Cartesian closed full subcategory within the category of continuous domains (resp., algebraic domains),   which  is essential for modeling the typed $\lambda$-calculus. In \cite{Lawson-Xu}, Lawson and Xu generalized  L-domains to sL-domains,  while Xu in \cite{xu} characterized sL-domains  by function spaces. Although Xu and Mao's study \cite{xu-mao-ab} delved into the representations of L-domains using abstract bases, they did not  mention the representation for sL-domains.   In this paper, different from the original concept in \cite{Jung1}, L-domains  need not  have a  least element. Following the work \cite{xu-mao-ab} and  \cite{Wu-Xu}, we will continue to study the connection between domain structures and CF-approximation spaces, aiming to extend Xu and Mao's representation approach from L-domains  to sL-domains within the framework of CF-approximation spaces.

The paper  is organized as follows: In Section 2, we recall some  basic notions and results in domain theory  that will be used later. In Section 3, we   introduce the notion of dense  CF-approximation subspaces, followed by an investigation into their properties concerning ordered structures. In Section 4, we introduce the concepts of  sL-approximation spaces and  L-approximation spaces, which are special types of CF-approximation spaces. Subsequently, we  prove that  every sL-approximation
space (resp.,  L-approximation space,  bc-approximation space) can induce  an sL-domain (resp.,  L-domain).
In Section 5,  we show that  every sL-domain (resp., L-domain, bc-domain) is order-isomorphic to the  collection of CF-closed sets of an sL-approximation space (resp., L-approximation
space, bc-approximation space) and  establish  equivalences  between relevant categories, such as the equivalence between category of sL-domains  and  that of  sL-approximation spaces.

\section{Preliminaries}

In this section,   we quickly recall some basic notions and results in Domain Theory. For  notions  not explicitly defined herein, the reader may refer to  \cite{Davey-Pri, Gierz,Jean-2013,jar2}. For some basic notions and results in category theory, please  refer to \cite{Well}.

For a set $U$ and $X\subseteq U$,  use $\mathcal {P} (U)$
to denote the power set of $U$, use $\mathcal{P}_{fin}(U)$ to denote the family of all  finite subsets of $U$ and  use $X^c$ to denote the complement
of $X$ in $U$. The symbol $F\subseteq_{fin} X$ means that  $F$ is a finite subset of $X$.

Let ($L$, $\leqslant$) be a poset. A \emph{principal ideal} (resp., \emph{principal filter}) of $L$
is a set of the form ${\da x}=\{y\in L\mid y\leqslant x\}$ (resp., ${\uparrow\! x}=\{y\in L\mid x\leqslant y\}$). For $A\subseteq L$, we write $\da A=\{y\in L\mid\exists \ x\in A, \,y\leqslant x\}$ and $\uparrow\! A=\{y\in L\mid\exists \ x\in A,\, x\leqslant
y\}$. A subset $A$ is a \emph{lower set} (resp., an \emph{upper set}) if $A=\da A$ (resp., $A={\uparrow}A$). We say that $z$ is a \emph{lower bound} (resp., an \emph{upper bound}) of $A$ if $A\subseteq \uparrow\!z$ (resp., $A\subseteq\downarrow\!z$).  A subset $B$ of $L$ is said to be {\em up-bounded} if $B$ has an
upper bound. The supremum of $A$ is the least upper bound of $A$, denoted by $\bigvee A$ or $\sup A$. The infimum of $A$ is the greatest lower bound of $A$, denoted by $\bigwedge A$ or $\inf A$. If $A\subseteq\downarrow\!x\subseteq L$, then we use $\bigvee_x A$ or  $\sup_x A$ to denote the supremum of $A$ in sub-poset $(\downarrow\!x,  \leqslant)$, where the order inherits from $L$. A nonempty subset $D$ of $L$ is \emph{directed} if every  nonempty finite subset of $D$ has an upper bound in $D$.
 A poset $L$ is a \emph{directed complete partially ordered set} ({\sl dcpo}, for short) if every directed subset of $L$ has a supremum.
%


Recall that in a poset $P$, we say that $x$ \emph{way-below}
$y$, written $x\ll y$,  if for any directed set $D$  having
a supremum with $\sup D\geqslant y$, there exists   some $d\in D$ such that $x\leqslant d$. If $x\ll x$, then $x$ is called a compact element of $P$. The set $\{x\in P\mid x\ll x\}$ is denoted by $K(P)$. 
The set $\{y\in P\mid  x\ll
y\}$ will be denoted $\uuar x$ and $\{y\in P\mid  y\ll x\}$ denoted
$\dda x$.  A poset $P$ is said to be
\emph{continuous} (resp., {\em algebraic}) if for all $ x\in P$, $\dda x$ is directed (resp., ${\downarrow} x\cap K(P)$ is directed) and  $x=\bigvee \dda x$ (resp., $x=\bigvee ({\downarrow} x\cap K(P))$). If a dcpo $P$ is  continuous (resp., algebraic), then $P$ is called a {\em continuous domain} (resp., {\em an algebraic domain}). A \emph{semilattice} (resp., {\em sup-semilattice}) is a poset in which every pair of elements has an infimum (resp., a supremum). A \emph{complete lattice} is a poset in which every subset has a supremum (equivalently, has an infimum). A subset $B$ of of a poset $P$ is called a {\em basis} of $P$ if for all $x\in P$, there exists   $B_x\subseteq B\cap\dda x$ such that
$B_x$ is directed and $\sup B_x=x$. It is well known that a poset $P$ is continuous if and only if it has a
 basis,  and that $P$ is algebraic  if and only if $K(P)$ is a
 basis.

Let $L$ and $P$ be dcpos, and $f: L \longrightarrow P$ a mapping. Then $f$ is called a {\em Scott continuous mapping} if for every directed subset $D\subseteq L$, $f(\bigvee D)=\bigvee f(D)$.

\begin{lm} {\em (\cite{Gierz})}\label{lm-way-rel} Let $P$ be a poset. Then for all $x, y, u, z\in P$,

{\em (1)} $x\ll y\Rightarrow x\leqslant y$;

{\em (2)} $u\leqslant x\ll y\leqslant z\Rightarrow u\ll z$;

{\em (3)} if $x\ll y$, $y\ll z$ and $x\vee y$ exists, then $x\vee y\ll z$; and

{\em (4)} if $P$ has the least element $\bot$, then $\bot\ll x$.

\end{lm}
\begin{lm}{\em (\cite{Gierz})}\label{pn-int}
If $P$ is a continuous poset, then the way-below relation
$\ll$ has
the  interpolation property:
   $x\ll z\Rightarrow \exists y\in P$ such that $x\ll y\ll
   z$.
\end{lm}
\begin{lm}\label{lm-bas-dom} {\rm (\cite[Lemma 4.1]{Wu-Xu})} Let $(L,\leqslant)$ be a continuous domain and $B\subseteq L$ a basis. If $(B, \leqslant)$ is a  sup-semilattice, then $(L, \leqslant)$ is a  sup-semilattic.
\end{lm}

\begin{lm}\label{lm-sup-st}{\rm (\cite[Lemma 2.3]{Wu-Xu-FIE})}Let $P$ be a poset, $A\subseteq P$ and $s, t$ be two upper bounds of $A$. If the supremums $\bigvee_{s}A$ and $\bigvee_{t}A$ exist and $\bigvee_{s}A\leqslant t$, then $\bigvee_{s}A=\bigvee_{t}A$
\end{lm}

\begin{dn} (\cite{Gierz,xu}) \label{dn-doms}
 (1)  A poset $P$ is called an {\em sL-cusl}, if for every $p\in P$, ${\downarrow} p$ is a sup-semilattice.

 (2)  A poset $P$ is called an {\em L-cusl}, if for every $p\in P$, ${\downarrow} p$ is a sup-semilattice with bottom element.

 (3)  A continuous domain $L$ is called a {\em continuous sL-domain} ({\em sL-domain}, for short), if for every $x\in L$, the set ${\downarrow} x$ is a  sup-semilattice in the order inherited from $L$.

 (4)  A continuous domain $L$ is called a  {\em continuous L-domain} ({\em L-domain}, for short), if for every $x\in L$, the set ${\downarrow} x$ is a complete lattice in the order inherited from $L$.

 (5)  A poset  $P$ is said {\em pointed}  if  $P$ has the least element.

 (6)  A poset  $P$ is called  a {\em cusl},  if every finite  up-bounded subset $A$ of $P$  has a  supremum.

 (7) A poset  $P$ is called  a {\em bc-poset}, if every up-bounded subset $B$ of  $P$  has a  supremum. If a bc-poset $P$ is also a continuous domain, then $P$ is called a {\em bc-domain}.
\end{dn}

 It is well-known that a bc-domain is a  pointed L-domain and a pointed sL-domain is a pointed L-domain. However, a pointed L-domain need not be a bc-domain.\\



In what follows,  we recall some contents on  binary relation, which are also the foundation of rough set theory.
\begin{dn} (\cite{gl-liu,Yang-Xu-Alg}) A
\emph{generalized approximation space}
(\emph{GA-space}, for short) is a pair $(U, R)$ consisting of a background set $U$ and
a   binary relation $R$ on $U$. For a GA-space $\ur$,
define $R_s, R_p: U\longrightarrow\mathcal {P}(U)$  by
\vskip3pt

\centerline{$R_s(x)=\{y\in U \mid   xRy\}$ $(\forall x\in U)$;\ \ $R_p(x)=\{y\in U \mid  yRx\}$ $(\forall x\in U).$}

\end{dn}

\begin{dn} (\cite{gl-liu,Yang-Xu-Alg})\label{dn-app}
 Let $(U,R)$ be a GA-space. Define $\underline{R},\overline{R}:\mathcal {P} (U)\longrightarrow
\mathcal {P}(U)$ by
\vskip3pt

\centerline {$\underline{R} (A)=\{x\in U\mid \ R_s(x)\subseteq A\}$;\ \ $\overline{R}(A)=\{x\in U\mid \ R_s(x)\cap A\neq \emptyset\}.$}
\vskip5pt

\noindent Operators $\underline{R}$ and $\overline{R}$  are called {\em the lower approximation operator} and {\em the upper
approximation operator} in $(U,R)$, respectively.
\end{dn}

 The notion of approximation operators is an  important concept in rough set theory. Approximation operators $\underline{R}$ and $\overline{R}$ are dual to each other; that is for very $A\subseteq U$,   $\overline{R}(A^c)=(\underline{R}(A))^c$. Each conclusion about lower approximation operators can be transformed into the counterpart of upper approximation operators. Since  upper
approximation operators play a key role in this paper,  we now list some useful lemmas of  upper
approximation operators as follows,  although they are easy to obtain by definition.
\begin{lm} {\em (\cite{gl-liu,Yang-Xu-Alg})}\label{lm-ula}
Let $(U,R)$ be a GA-space.

{\em(1)} For every $\{A_{i}\mid i\in I\}\subseteq \mathcal {P}(U)$,
$\overline{R}(\bigcup_{i\in
I}A_{i})=\bigcup_{i\in I}\overline{R}(A_{i}).$

{\em(2)} For all $x\in U, \ \overline{R}(\{x\})=R_p(x)$.
\end{lm}

 A binary relation $R\subseteq U\times U$ on a set $U$ is said to be {\em reflexive} if for all $x\in U$, $x R x$; and
$R$ is said to be \emph{transitive} if $xRy$ and $yRz$ imply $xRz$ for all
$x,y,z\in U$. A binary relation $R$ is called a {\em preorder} if it is  reflexive and transitive. For a binary relation $R\subseteq U\times U$ and $ V\subseteq U$, we use $R|_{V}$ to denote $R\cap (V\times V)$.
\begin{lm}\label{lm-ref-lm2-tr-up}{\em (\cite{Wu-Xu})}
Let $(U,R)$ be a GA-space.

{\em(1)} $R$ is reflexive iff   for
all $X\subseteq U$,
 $X\subseteq \overline{R}(X)$.

{\em(2)} $R$ is transitive iff for
all $X\subseteq U$,
 $\overline{R}(\overline{R}(X))\subseteq \overline{R}(X)$.

{\em(3)} If $R$ is  transitive and $A, B\subseteq U$, then $\overline{R}(B)\subseteq\overline{R}(A)$
whenever $B\subseteq\overline{R}(A)$.

{\em(4)}  If $R$ is a  preorder, then the operator $\overline{R}$ is a topological closure  operator on $U$.
\end{lm}

\section{Dense  CF-approximation subspaces}

In \cite{Wu-Xu}, Wu and Xu introduced CF-approximation spaces and gave representations of continuous domains by them. In this section we first recall some key notions and results  in \cite{Wu-Xu}. Moreover, we will introduce the notion of dense  CF-approximation subspaces and study the ordered structures of  their CF-closed sets.

\begin{dn}\label{dn-cf-ga}(\cite{Wu-Xu}) Let $(U, R)$ be a GA-space, let $R$ be a transitive relation and let $\emptyset\ne\mathcal{F}\subseteq \mathcal{P}_{fin}(U)$. Then $(U, R, \mathcal{F})$ is called a {\em generalized approximation space with a consistent family of finite subsets} ({\em  CF-approximation space}, for short) if for every $ F\in \mathcal{F}$,  whenever $K\subseteq_{fin} \overline{R}(F)$, there always exists $G\in \mathcal{F}$ such that $K\subseteq \overline{R}(G)$ and $G\subseteq \overline{R}(F)$.
\end{dn}
\begin{dn}\label{dn-cf-cl}  (\cite{Wu-Xu})  Let $(U, R, \mathcal{F})$ be a CF-approximation space and $E\subseteq U$. Then $E$ is called {\em a CF-closed set} of $(U, R, \mathcal{F})$ if for all $ K\subseteq_{fin}E$, there always exists $F\in \mathcal{F}$ such that $K\subseteq \overline{R}(F)\subseteq E$ and $F\subseteq E$. The collection of all CF-closed sets of $(U, R, \mathcal{F})$ is denoted by $\mathfrak{C}(U, R, \mathcal{F})$.
\end{dn}

\begin{lm}{\em (\cite[Proposition 3.7]{Wu-Xu})}\label{pn2-cf-clo} For  a CF-approximation space $(U, R, \mathcal{F})$, the following statements hold:

$(1)$ For every $F\in \mathcal{F}$, $\overline{R}(F)\in \mathfrak{C}(U, R, \mathcal{F})$.

$(2)$ If $E\in \mathfrak{C}(U, R, \mathcal{F})$, $A\subseteq E$, then $\overline{R}(A)\subseteq E$.

$(3)$ If $\{E_{i}\}_{i\in I}\subseteq\mathfrak{C}(U, R, \mathcal{F})$ is a directed family, then $\bigcup_{i\in I}E_{i}\in \mathfrak{C}(U, R, \mathcal{F})$.

$(4)$ $(\mathfrak{C}(U, R, \mathcal{F}), \subseteq)$ is a dcpo.
\end{lm}

\begin{lm}{\em (\cite[Proposition 3.8]{Wu-Xu})}\label{pn3-cf-cl}
For  a CF-approximation space $(U, R, \mathcal{F})$ and $E\subseteq U$, the following statements are equivalent:

$(1)$ $E\in \mathfrak{C}(U, R, \mathcal{F})$;

$(2)$ The family $\mathcal{A}=\{\overline{R}(F)\mid F\in \mathcal{F},  F\subseteq E\}$ is directed and $E=\bigcup\mathcal{A}$;

$(3)$ There exists a family $\{F_{i}\}_{i\in I}\subseteq \mathcal{F}$ such that $\{\overline{R}(F_{i})\}_{i\in I}$ is directed, and $E=\bigcup_{i\in I}\overline{R}(F_{i})$;

$(4)$  There always  exists $F\in \mathcal{F}$ such that $K\subseteq \overline{R}(F)\subseteq E$ whenever $ K\subseteq_{fin}E$.
\end{lm}

\begin{ex}  (\cite[Theorem 3.13]{Wu-Xu})\label{tm-CDO-CFGA} Let $(L, \leqslant)$ be a continuous domain.  Write $R_{L}$  for the way-below relation ``$\ll$" of $(L, \leqslant)$ and write $\mathcal{F}_{L}=\{F\subseteq_{fin} L\mid F\  \mbox{has a top element } c_F\}$. Then

$(1)$ $(L, R_{L},  \mathcal{F}_{L})$ is a CF-approximation space;

$(2)$ $\mathfrak{C}(L, R_{L},  \mathcal{F}_{L})=\{\dda x\mid x\in L\}$.
\end{ex}

 We call $(L, R_{L},  \mathcal{F}_{L})$  the {\em induced CF-approximation space} of the continuous domain $L$.

\begin{lm}{\em (\cite[Theorem 3.9, Corollary 3.10]{Wu-Xu})}\label{tm-cfga-way}
Let $E_{1}, E_{2}\in \mathfrak{C}(U, R, \mathcal{F})$.  Then $ E_{1}\ll E_{2}$ if and only if there exists $F\in \mathcal{F}$ such that $ E_{1}\subseteq \overline{R}(F)$ and  $F\subseteq E_{2}$;

%
\end{lm}

\begin{cy}{\em (\cite[Theorem 3.11]{Wu-Xu})}\label{pn-bas-clo} Let $(U, R, \mathcal{F})$ be a CF-approximation space. Then $\{\overline{R}(F)\mid F\in \mathcal{F}\}$ is a basis of $(\mathfrak{C}(U, R, \mathcal{F}), \subseteq)$; and thus $(\mathfrak{C}(U, R, \mathcal{F}), \subseteq)$ is a continuous domain.
\end{cy}

By Lemma \ref{lm-ref-lm2-tr-up},   if $(U, R)$ is a GA-space with $R$ being a preorder and $\emptyset\ne\mathcal{F}\subseteq \mathcal{P}_{fin}(U)$, then $\overline{R}$ is a topological closure operator and $(U, R, \mathcal{F})$ is a CF-approximation space. So we have the following definition.

\begin{dn}\label{dn-cf-tga} (\cite{Wu-Xu})  Let $(U, R)$ be a GA-space and $\emptyset\ne\mathcal{F}\subseteq \mathcal{P}_{fin}(U)$.  If $R$ is a preorder, then the CF-approximation space $(U, R, \mathcal{F})$ is said to be {\em  topological}.
\end{dn}

\begin{lm}\label{lm-top-alg}  Let $(U, R, \mathcal{F})$ be a  topological CF-approximation space. Then $\{\overline{R}(F)\mid F\in \mathcal{F}\}=K(\mathfrak{C}(U, R, \mathcal{F}))$; thus  $(\mathfrak{C}(U, R, \mathcal{F}), \subseteq)$ is an algebraic domain.
\end{lm}

 \begin{ex}\label{ex-alg-top-cf-apps} (\cite{Wu-Xu})
 Let $(L, \leqslant)$ be an algebraic  domain. Write $R_{K(L)}=\leqslant_{K(L)}$ and write $\mathcal{F}_{K(L)}=\{F\subseteq_{fin} K(L)\mid F\  \mbox{has a top element } c_F\}$, where $\leqslant_{K(L)}$ is the restriction of $\leqslant$ to $K(L)$. Then

$(1)$ $(K(L), R_{K(L)},  \mathcal{F}_{K(L)})$ is a topological CF-approximation space;

$(2)$ $\mathfrak{C}(K(L), R_{K(L)},  \mathcal{F}_{K(L)})=\{{\downarrow} x\cap K(L)\mid x\in L\}$.
 \end{ex}

We call $(K(L),  R_{K(L)}, \mathcal{F}_{K(L)})$ the {\em induced topological CF-approximation space} of the algebraic domain $L$.

In what follows, we introduce the notions of   CF-approximation subspaces  and dense CF-approximation subspaces. We will  study the connection between    CF-approximation spaces and their   dense subspaces  in terms of ordered structure.
\begin{dn}
Let $(U, R, \mathcal{F})$ be a CF-approximation space, $V\subseteq U$ and $\mathcal{G}\subseteq \mathcal{F}\cap \mathcal{P}(V)$. Write $R|_V=R\cap (V\times V)$. Then $(V, R|_V, \mathcal{G})$ is called  a {\em CF-approximation subspace} of  $(U, R, \mathcal{F})$ if $(V, R|_V, \mathcal{G})$ itself is a CF-approximation space.
\end{dn}
\begin{lm}\label{lm-dense-cf} Let $(U, R, \mathcal{F})$ be a CF-approximation space, $V\subseteq U$ and $\mathcal{G}\subseteq \mathcal{F}\cap \mathcal{P}(V)$. If for every $K\subseteq_{fin}U$ and $F\in \mathcal{F}$ with $K\subseteq \overline{R}(F)$, there exists $G\in \mathcal{G}$ such that $K\subseteq \overline{R}(G)$ and $G\subseteq\overline{R}(F)$. Then $(V, R|_V, \mathcal{G})$ is a  CF-approximation subspace of $(U, R, \mathcal{F})$, where $R|_V=R\cap(V\times V)$.
\end{lm}
\noindent {\bf Proof.}
It suffices to prove that $(V, R|_V, \mathcal{G})$ is a CF-approximation space. Clearly, $R|_V=R\cap(V\times V)$ is a transitive relation on $V$. Let $K\subseteq_{fin} V$ and $G\in \mathcal{G}$. If $K\subseteq \overline{R|_V}(G)=\overline{R}(G)\cap V$, then $K\subseteq_{fin} \overline{R}(G)$. It follows from the assumption that there exists $H\in \mathcal{G}$ such that $K\subseteq\overline{R}(H)$ and $H\subseteq\overline{R}(G)$. Noticing that $K, H\in \mathcal{P}(V)$ , we have

\centerline{$K\subseteq\overline{R}(H)\cap V=\overline{R|_V}(H)$, $H\subseteq\overline{R}(G)\cap V=\overline{R|_V}(G)$.}
\noindent Thus, $(V, R|_V, \mathcal{G})$ is a  CF-approximation space.\hfill$\Box$\\

In view of Lemma \ref{lm-dense-cf}, we give the following definition.
\begin{dn}\label{dn-dence-sub}
Let $(U, R, \mathcal{F})$ be a CF-approximation space, $V\subseteq U$ and $\mathcal{G}\subseteq \mathcal{F}\cap \mathcal{P}(V)$. If for every $K\subseteq_{fin}U$ and $F\in \mathcal{F}$ with $K\subseteq \overline{R}(F)$, there exists $G\in \mathcal{G}$ such that $K\subseteq \overline{R}(G)$ and $G\subseteq\overline{R}(F)$, then we say the CF-approximation subspace $(V, R|_V, \mathcal{G})$ is {\em dense} in $(U, R, \mathcal{F})$, or, a {\em dense subspace} of $(U, R, \mathcal{F})$.
\end{dn}

\begin{pn}\label{pn-dens-clo}
Let $(V, R|_V, \mathcal{G})$ be a dense CF-approximation subspace of  $(U, R, \mathcal{F})$. Then
$$(\mathfrak{C}((U, R, \mathcal{F})), \subseteq)\cong(\mathfrak{C}(V, R|_V, \mathcal{G}), \subseteq).$$
\end{pn}
\noindent{\bf Proof.}   we divide the  proof into several steps.

Step 1. Assert that for all $E\in \mathfrak{C}(U, R, \mathcal{F})$, we have $E\cap V\in \mathfrak{C}(V, R|_V, \mathcal{G})$.

Let $E\in \mathfrak{C}(U, R, \mathcal{F})$. If $K\subseteq_{fin} E\cap V$, then $K\subseteq_{fin} E$. It follows from $E\in \mathfrak{C}(U, R, \mathcal{F})$ that there exists $F\in \mathcal{F}$ such that $K\subseteq \overline{R}(F)\subseteq E$ and $F\subseteq E$. By definition \ref{dn-cf-cl}, there exists $G\in \mathcal{G}$ such that $ K\subseteq \overline{R}(G)$ and $G\subseteq \overline{R}(F)\subseteq E$. Noticing that $K, G\in \mathcal{P}(V)$ , we have

\centerline{$K\subseteq\overline{R}(G)\cap V=\overline{R|_V}(G)\subseteq E\cap V$, $G\subseteq E\cap V$.}
\noindent Thus, $E\cap V\in \mathfrak{C}(V, R|_V, \mathcal{G})$.

Step 2. Assert that for $E_1, E_2\in \mathfrak{C}(U, R, \mathcal{F})$,  $E_1\cap V\subseteq E_2\cap V$ implies $E_1\subseteq E_2$.

 Let $x\in E_1$. Then there exists $F\in \mathcal{F}$ such that $x\in \overline{R}(F)\subseteq E_1$ and $F\subseteq E_1$. Since $(V, R|_V, \mathcal{G})$ is a dense CF-approximation subspace,  there exists $G\in \mathcal{G}$ such that $x\in \overline{R}(G)$ and $G\subseteq\overline{R}(F)\subseteq E_1$. Noticing that $G\subseteq V$, we have

 \centerline{$G\subseteq E_1\cap V\subseteq E_2\cap V\subseteq E_2$.}
 \noindent Thus $x\in \overline{R}(G)\subseteq E_2$. By arbitrariness of $x\in E_1$, we have $E_1\subseteq E_2$.

Step 3. Assert that  $\mathfrak{C}(V, R|_V, \mathcal{G})=\{E\cap V\mid E\in\mathfrak{C}(U, R, \mathcal{F})\}$.

 Let $S\in \mathfrak{C}(V, R|_V, \mathcal{G})$. Then by Proposition \ref{pn3-cf-cl}(2), we have that
\begin{align*}
S&=\bigcup\{\overline{R|_V}(G)\mid G\in \mathcal{G}, G\subseteq S\}\\
&=\bigcup\{\overline{R}(G)\cap V\mid G\in \mathcal{G}, G\subseteq S\}\\
&=(\bigcup\{\overline{R}(G)\mid G\in \mathcal{G}, G\subseteq S\})\cap V
\end{align*}
\noindent and that $(\{\overline{R}(G)\cap V\mid G\in \mathcal{G}, G\subseteq S\}, \subseteq)$ is directed. By Step 2, we have that $(\{\overline{R}(G)\mid G\in \mathcal{G}, G\subseteq S\}, \subseteq)$ is directed. Thus, $\bigcup\{\overline{R}(G)\mid G\in \mathcal{G}, G\subseteq S\}\in \mathfrak{C}(U, R, \mathcal{F})$. By arbitrariness of $S\in \mathfrak{C}(V, R|_V, \mathcal{G})$, we have $\mathfrak{C}(V, R|_V, \mathcal{G})\subseteq\{E\cap V\mid E\in\mathfrak{C}(U, R, \mathcal{F})\}$. It directly follows from Step 1 that $\{E\cap V\mid E\in\mathfrak{C}(U, R, \mathcal{F})\}\subseteq\mathfrak{C}(V, R|_V, \mathcal{G})$. Thus,
\vskip 3pt
\centerline{$\mathfrak{C}(V, R|_V, \mathcal{G})=\{E\cap V\mid E\in\mathfrak{C}(U, R, \mathcal{F})\}$.}
\vskip 3pt
Step 4. Assert that $(\mathfrak{C}((U, R, \mathcal{F})), \subseteq)\cong(\mathfrak{C}(V, R|_V, \mathcal{G}), \subseteq)$.

 Define  a mapping $f: (\mathfrak{C}((U, R, \mathcal{F})), \subseteq)\longrightarrow(\mathfrak{C}(V, R|_V, \mathcal{G}), \subseteq)$ by
$$\forall E\in \mathfrak{C}((U, R, \mathcal{F})), f(E)=E\cap V\in \mathfrak{C}(V, R|_V, \mathcal{G}).$$
By Step 2 and Step 3, we have that $f$ is  a bijection, $f$ and $f^{-1}$ are order-preserving.  So,  $(\mathfrak{C}((U, R, \mathcal{F})), \subseteq)\cong(\mathfrak{C}(V, R|_V, \mathcal{G}), \subseteq)$.\hfill$\Box$\\

It is easy to see that  the topological CF-approximation space $(K(L),  R_{K(L)}, \mathcal{F}_{K(L)})$ is a dense subspace of  $(L, R_{L},  \mathcal{F}_{L})$  for every algebraic domain  $(L, \leqslant)$. By Lemma \ref{lm-top-alg}, Example \ref{ex-alg-top-cf-apps} and Proposition \ref{pn-dens-clo}, we get a new representation for algebraic domains as follows.
\begin{tm}{\em (Representation Theorem for algebraic domains)}
A poset $(L, \leqslant)$ is an algebraic domain if and only if there exists a CF-approximation $(U, R, \mathcal{F})$  with a topological CF-approximation space as a dense subspace such that $(\mathfrak{C}(U, R, \mathcal{F}), \subseteq)\cong (L, \leqslant)$.
\end{tm}

\section{sL-approximation spaces}

In this section, we will introduce  two   suitable  types of  CF-approximation spaces that can generate sL-domains and L-domains, respectively.
It is worth noting that a continuous domain is an sL-domain if its principal ideals are sup-semilattices. Therefore,  we first  consider a   CF-approximation subspace induced by a CF-closed set of the original CF-approximation space, where the family of  CF-closed sets of the subspace is exactly a principal ideal of that  of the original CF-approximation space.

\begin{pn}\label{pn-cons-idl} Let $(U, R, \mathcal{F})$ be a CF-approximation space, $E\in \mathfrak{C}(U, R, \mathcal{F})$ and $\mathcal{F}_{E}=\{F\in \mathcal{F}\mid F\subseteq E\}$. Then

$(1)$ $(E,  R|_E, \mathcal{F}_{E})$ is  a CF-approximation space; and

$(2)$ $\mathfrak{C}(E,  R|_E, \mathcal{F}_{E})=\downarrow\! E=\{H\in \mathfrak{C}(U, R, \mathcal{F})\mid H\subseteq E\}$.
\end{pn}
\noindent{\bf Proof.} (1) Let $F\in \mathcal{F}_{E}$ and $K\subseteq_{fin}E$. If $K\subseteq \overline{R|_E}(F)=\overline{R}(F)\cap E$, then $K\subseteq_{fin}\overline{R}(F)$. By Definition \ref{dn-cf-cl}, there exists $G\in F$ such that $K\subseteq \overline{R}(G)$ and $G\subseteq \overline{R}(F)$. It follows from $F\in \mathcal{F}_{E}$  and Lemma \ref{pn2-cf-clo}(2) that $F\subseteq E$ and $G\subseteq\overline{R}(F)\subseteq E$. Thus $G\in \mathcal{F}_{E}$. Noticing that $K\subseteq E$ and $G\subseteq E$, we have $K\subseteq\overline{R}(G)\cap E=\overline{R|_E}(G)$ and $G\subseteq\overline{R}(F)\cap E=\overline{R|_E}(F)$. Thus $(E,  R|_E, \mathcal{F}_{E})$ is a CF-approximation space.

(2) For all $F\in \mathcal{F}_{E}$, it follows from Lemma \ref{pn2-cf-clo}(2) that $\overline{R}(F)\subseteq E$. Thus $\overline{R|_E}(F)=\overline{R}(F)\cap E=\overline{R}(F)$. Let $H\in \mathfrak{C}(E,  R|_E, \mathcal{F}_{E})$. Clearly, $H\subseteq E$. If $K\subseteq_{fin} H$, then there exists $F\in \mathcal{F}_{E}$ such that $K\subseteq\overline{R|_E}(F)=\overline{R}(F)\subseteq H$ and $F\subseteq H$, showing that  $H\in \mathfrak{C}(U, R, \mathcal{F})$. Thus  $\mathfrak{C}(E,  R|_E, \mathcal{F}_{E})\subseteq{\downarrow} E$.  Conversely, let $H\in \mathfrak{C}(U, R, \mathcal{F})$ and $ H\subseteq E$. If $K\subseteq_{fin}H$, then there exists $F\in \mathcal{F}$ such that $K\subseteq \overline{R}(F)\subseteq H$ and $F\subseteq H$. It follows from $H\subseteq E$  that $F\in \mathcal{F}_E$ and  $\overline{R|_E}(F)=\overline{R}(F)$.  Therefore $K\subseteq \overline{R|_E}(F)\subseteq H$. Thus $H\in \mathfrak{C}(E,  R|_E, \mathcal{F}_{E})$, showing that $\mathfrak{C}(E,  R|_E, \mathcal{F}_{E})\supseteq\downarrow\! E$. Thus

\vskip 3pt
\centerline{$\mathfrak{C}(E,  R|_E, \mathcal{F}_{E})=\downarrow\! E=\{H\in \mathfrak{C}(U, R, \mathcal{F})\mid H\subseteq E\}$.}
\hfill$\Box$

The CF-approximation space $(E,  R|_E, \mathcal{F}_{E})$ is   called the {\em induced principal CF-approximation subspace} of the CF-closed set $E$.
 By Proposition \ref{pn-cons-idl}, $\mathfrak{C}(E,  R|_E, \mathcal{F}_{E})=\downarrow\! E$, exactly a principal ideal of $(\mathfrak{C}(U, R, \mathcal{F}), \subseteq)$.

 \begin{dn} {(\cite{Rong-Xu})}
 Let $U$ be a  nonempty set, $\mathcal{F}\neq\emptyset$ be a family of subsets of $U$.
We use $\mathcal{F}^\sharp$  denote the set of all the nonempty finite unions of $\mathcal{F}$, that is
\vskip 3pt
 \centerline{$\mathcal{F}^\sharp=\{A\mid \exists \mathcal{A}\subseteq_{fin} \mathcal{F}, \mathcal{A}\neq\emptyset \mbox{ such that } A=\cup\mathcal{A}\}.$}
\noindent Then $\mathcal{F}^\sharp$ is called the {\em join saturation} of $\mathcal{F}$.
\end{dn}
\begin{dn} Let $(U, R, \mathcal{F})$ be a CF-approximation space, $K\in \mathcal{F}^{\sharp}\cup\{\emptyset\}$, $F\in \mathcal{F}$ with $\overline{R}(K)\subseteq \overline{R}(F)$. Define $\mathcal{S}(K, F)$ to be the collection of all  $G\in \mathcal{F}$ satisfies the following conditions:

{(sL-1)} $\overline{R}(K)\subseteq\overline{R}(G)\subseteq\overline{R}(F)$;

{ (sL-2)} $\forall G^{\prime} \in \mathcal{F}$, $\overline{R}(K)\subseteq\overline{R}(G^{\prime})\subseteq\overline{R}(F)\Longrightarrow\overline{R}(G)\subseteq\overline{R}(G^{\prime})$.

\noindent Define $\mathcal{S}^{*}(K, F)=\{G\in \mathcal{S}(K, F)\mid G\subseteq \overline{R}(F)\}$.
\end{dn}
\begin{rk} \label{rk-*subseteq} For a family $\mathcal{F}$, $\emptyset \in \mathcal{F}\Leftrightarrow \emptyset \in \mathcal{F}^{\sharp}\Leftrightarrow \mathcal{F}^{\sharp}\cup\{\emptyset\}=\mathcal{F}^{\sharp}$. For a CF-approximation space $(U, R, \mathcal{F})$, if $K\in \mathcal{F}^{\sharp}\cup\{\emptyset\}$ and $F\in \mathcal{F}$ with $\overline{R}(K)\subseteq \overline{R}(F)$, then $\mathcal{S}^{*}(K, F)\subseteq\mathcal{S}(K, F)$.
\end{rk}

\begin{pn}\label{1-pn-M-sup}
Let $(U, R, \mathcal{F})$ be a CF-approximation space, $F, F_1, F_2, F_3\in \mathcal{F}$, $M, M_{1}, M_{2}\in \mathcal{F}^{\sharp}\cup\{\emptyset\}$. Then the following statements hold.

$(1)$ If     $\overline{R}(M)\subseteq\overline{R}(F)$, $G_{1}\in \mathcal{S}(M, F)$ and $G_{2}\in \mathcal{F}$, then $G_{2}\in \mathcal{S}(M, F)$ iff $\overline{R}(G_{1})=\overline{R}(G_{2})$.

$(2)$ If   $\overline{R}(M)\subseteq\overline{R}(F_{1})\subseteq\overline{R}(F_{2})$ and $\mathcal{S}(M, F_{1})\ne \emptyset\ne\mathcal{S}(M, F_{2})$, then $\mathcal{S}(M, F_{1})=\mathcal{S}(M, F_{2})$.

$(3)$  If    $\overline{R}(M_{1})\subseteq\overline{R}(M_{2})\subseteq\overline{R}(F)$ and $G_{i}\in
\mathcal{S}(M_{i}, F)\ (i=1, 2)$, then $\overline{R}(G_{1})\subseteq\overline{R}(G_{2})$.

$(4)$ If   $\overline{R}(M)\subseteq\overline{R}(F_{1})\cap\overline{R}(F_{2})$, $\overline{R}(F_{1}\cup F_{2})\subseteq \overline{R}(F_{3})$ and $\mathcal{S}(M, F_{i})\neq\emptyset\ (i=1, 2, 3)$, then $\mathcal{S}(M, F_{1})=\mathcal{S}(M, F_{2})$.

\end{pn}

\noindent{\bf Proof.}
(1) If $G_{2}\in \mathcal{S}(M, F)$,  then by (sL-1), we have  $\overline{R}(M)\subseteq\overline{R}(G_{i})\subseteq\overline{R}(F)$ $(i=1, 2)$. It follows from (sL-2) that $\overline{R}(G_{1})\subseteq\overline{R}(G_{2})$ and $\overline{R}(G_{2})\subseteq\overline{R}(G_{1})$. Thus $\overline{R}(G_{1})=\overline{R}(G_{2})$. Conversely, if $\overline{R}(G_{2})=\overline{R}(G_{1})$, then it follows from $G_{1}\in \mathcal{S}(M, F)$ that $G_{2}$ satisfies (sL-1) and (sL-2), showing that $G_{2}\in \mathcal{S}(M, F)$.

(2) Let $G_{1}\in \mathcal{S}(M, F_{1})$ and $G_{2}\in \mathcal{S}(M, F_{2})$. Then  $\overline{R}(M)\subseteq\overline{R}(G_{1})\subseteq\overline{R}(F_1)\subseteq\overline{R}(F_{2})$ and $\overline{R}(M)\subseteq\overline{R}(G_{2})\subseteq\overline{R}(F_{2})$.  By  (sL-2) and $G_{2}\in \mathcal{S}(M, F_{2})$, we have $\overline{R}(G_{2})\subseteq\overline{R}(G_{1})\subseteq\overline{R}(F_1)$. It follows from   $\overline{R}(M)\subseteq\overline{R}(G_{2})\subseteq\overline{R}(F_1)$,  $G_{1}\in \mathcal{S}(M, F_{1})$ and (sL-2) that $\overline{R}(G_{1})\subseteq\overline{R}(G_{2})$. Thus $\overline{R}(G_{1})=\overline{R}(G_{2})$. By (1), we have  $\mathcal{S}(M, F_{1})=\mathcal{S}(M, F_{2})$.

(3) It follows from   $G_{2}\in \mathcal{S}(M_{2}, F)$ that $\overline{R}(M_{1})\subseteq\overline{R}(M_{2})\subseteq\overline{R}(G_{2})\subseteq\overline{R}(F)$. By $G_{1}\in
\mathcal{S}(M_{1}, F)$ and (sL-2), we have  $\overline{R}(G_{1})\subseteq\overline{R}(G_{2})$.

(4) It follows from $\overline{R}(F_{1}\cup F_{2})\subseteq \overline{R}(F_{3})$, $\mathcal{S}(M, F_{i})\neq\emptyset\ (i=1, 2, 3)$ and (2) that $\mathcal{S}(M, F_{1})=\mathcal{S}(M, F_{3})$ and $\mathcal{S}(M, F_{2})=\mathcal{S}(M, F_{3})$. Thus, $\mathcal{S}(M, F_{1})=\mathcal{S}(M, F_{2})$.
\hfill$\Box$

\begin{pn} \label{pn-sup-semi}Let $(U, R, \mathcal{F})$ be a CF-approximation space, $E\in \mathfrak{C}(U, R, \mathcal{F})$. If $\mathcal{S}^{*}(K, F)\neq\emptyset$ for all $K\in \mathcal{F}^{\sharp}$ and  $F\in \mathcal{F}$ with $\overline{R}(K)\subseteq \overline{R}(F)$,  then $(\mathfrak{C}(E,  R|_E, \mathcal{F}_{E}), \subseteq)$ is a sup-semilattice.
\end{pn}
\noindent{\bf Proof.} It follows from Proposition \ref{pn-cons-idl} that $\mathfrak{C}(E,  R|_E, \mathcal{F}_{E})=\downarrow\! E$ and $\{\overline{R}(F)\mid F\in \mathcal{F}, F\subseteq E\}$ is the basis of $(\downarrow\! E, \subseteq)$. To prove $(\downarrow\! E, \subseteq)$ is a sup-semilattice, by Lemma \ref{lm-bas-dom}, it suffices to prove $(\{\overline{R}(F)\mid F\in \mathcal{F}, F\subseteq E\}, \subseteq)$ is a sup-semilattice. Let $F_1, F_2\in \mathcal{F}_{E}$. Then $F_1\cup F_2\subseteq_{fin}E$. By Definition \ref{dn-cf-cl}, there exists $F_{3}\in \mathcal{F}$ such that $F_1\cup F_2\subseteq \overline{R}(F_{3})\subseteq E$ and $F_{3}\subseteq E$. Let  $G\in \mathcal{S}^{*}(F_1\cup F_2, F_3)$. Then $\overline{R}(F_1)\cup \overline{R}(F_2)=\overline{R}(F_1\cup F_2)\subseteq\overline{R}(G)
\subseteq\overline{R}(F_{3})$ and $G\subseteq \overline{R}(F_{3})\subseteq E$, showing that $\overline{R}(G)$ is a  upper bound of $\overline{R}(F_1)$ and $ \overline{R}(F_2)$ in $(\{\overline{R}(F)\mid F\in \mathcal{F}, F\subseteq E\}, \subseteq)$. Let $F_{4}\in \mathcal{F}_E$ and $\overline{R}(F_1)\cup \overline{R}(F_2)\subseteq\overline{R}(F_4)$. It follows from $F_3\cup F_4\subseteq_{fin} E$ and Definition \ref{dn-cf-cl} that there exists $F^{\prime}\in \mathcal{F}$ such that $F_3\cup F_4\subseteq \overline{R}(F^{\prime})$ and $F^{\prime}\subseteq E$. Thus  $\overline{R}(F_{3}\cup F_{4})\subseteq \overline{R}(F^{\prime})$.   By Remark \ref{rk-*subseteq} and Proposition \ref{1-pn-M-sup}(4),  we have $\mathcal{S}(F_1\cup F_2, F_4)=\mathcal{S}(F_1\cup F_2, F_3)$. For $G^{\prime}\in \mathcal{S}^{*}(F_1\cup F_2, F_4)$, by Proposition \ref{1-pn-M-sup}(1),  $\overline{R}(G)=\overline{R}(G^{\prime})\subseteq\overline{R}(F_4)$,  showing that $\overline{R}(G)$ is the least upper bound of $\overline{R}(F_1)$ and $ \overline{R}(F_2)$ in $(\{\overline{R}(F)\mid F\in \mathcal{F}, F\subseteq E\}, \subseteq)$. Thus $\downarrow\! E=\{H\in \mathfrak{C}(U, R, \mathcal{F})\mid H\subseteq E\}$ is a sup-semilattice.\hfill$\Box$

\begin{cy} \label{cy-sup-sL} Let $(U, R, \mathcal{F})$ be a CF-approximation space. If $\mathcal{S}^{*}(K, F)\neq\emptyset$ for all $K\in \mathcal{F}^{\sharp}$ and  $F\in \mathcal{F}$ with $\overline{R}(K)\subseteq \overline{R}(F)$,  then $(\mathfrak{C}(U, R, \mathcal{F}), \subseteq) $ is an sL-domain.
\end{cy}
\noindent{\bf Proof.}  By Corollary \ref{pn-bas-clo}, we have that  $(\mathfrak{C}(U, R, \mathcal{F}), \subseteq)$ is a continuous domain. By Proposition \ref{pn-cons-idl} and \ref{pn-sup-semi},
we have that $(\mathfrak{C}(U, R, \mathcal{F}), \subseteq) $ is an sL-domain.\hfill$\Box$\\

Based on  Corollary \ref{cy-sup-sL}, we give the definition of  ultra sL-approximation spaces.

\begin{dn}\label{dn-sup-sL-App}
Let $(U, R, \mathcal{F})$ be a CF-approximation space, $\mathcal{F}^{\sharp}$ the join saturation of $\mathcal{F}$. If $\mathcal{S}^{*}(M, F)\neq\emptyset$ for all $M\in \mathcal{F}^{\sharp}$ and  $F\in \mathcal{F}$ with $\overline{R}(M)\subseteq \overline{R}(F)$;   that is to say, there exists $G\in \mathcal{F}$ satisfies the following conditions  for all $M\in \mathcal{F}^{\sharp}$ and  $F\in \mathcal{F}$ with $\overline{R}(M)\subseteq \overline{R}(F)$:

{(sL-$1^{\prime})$} $\overline{R}(M)\subseteq\overline{R}(G)$, $G\subseteq\overline{R}(F)$; and

{ (sL-2)} $\forall G^{\prime} \in \mathcal{F}$, $\overline{R}(M)\subseteq\overline{R}(G^{\prime})\subseteq\overline{R}(F)\Longrightarrow\overline{R}(G)\subseteq\overline{R}(G^{\prime})$,

\noindent then   $(U, R, \mathcal{F})$ is called an {\em ultra  sL-approximation space}.
\end{dn}
By Corollary \ref{cy-sup-sL},  for an ultra sL-approximation space $(U, R, \mathcal{F})$, one has  that  $(\mathfrak{C}(U, R, \mathcal{F}), \subseteq)$ is an sL-domain.\\

The conditions of  ultra sL-approximation spaces are too strong to be satisfied by all induced CF-approximation spaces of sL-domains (see Example \ref{ex-sl<ultra} in next section).  So,  we consider   weaker conditions of ultra sL-approximation spaces
which  can  also guarantee the collection of CF-closed sets forms an sL-domain.\\

By shifting the precondition ``$\overline{R}(M)\subseteq\overline{R}(F)$" in Definition \ref{dn-sup-sL-App} to ``$M\subseteq\overline{R}(F)$" and shifting the condition (sL-$1^{\prime}$) in Definition \ref{dn-sup-sL-App} to (sL-$1$), we give the following definition.
\begin{dn}\label{dn-sL-App}
Let $(U, R, \mathcal{F})$ be a CF-approximation space.  If $\mathcal{S}(M, F)\neq\emptyset$ for all $M\in \mathcal{F}^{\sharp}$ and  $F\in \mathcal{F}$ with  $M\subseteq\overline{R}(F)$,   that is  there exists $G\in \mathcal{F}$ satisfies the following conditions for all $M\in \mathcal{F}^{\sharp}$ and  $F\in \mathcal{F}$ with  $M\subseteq\overline{R}(F)$:

{(sL-$1)$} $\overline{R}(M)\subseteq\overline{R}(G)\subseteq\overline{R}(F)$; and

{(sL-2)} $\forall G^{\prime} \in \mathcal{F}$, $\overline{R}(M)\subseteq\overline{R}(G^{\prime})\subseteq\overline{R}(F)\Longrightarrow\overline{R}(G)\subseteq\overline{R}(G^{\prime})$,

\noindent then  $(U, R, \mathcal{F})$ is called an {\em sL-approximation space}.
\end{dn}

It is clear that (sL-$1^{\prime}$) implies (sL-1), that $M\subseteq\overline{R}(F)$ implies $\overline{R}(M)\subseteq\overline{R}(F)$  and that $\mathcal{S}^*(M, F)\neq\emptyset$ implies $\mathcal{S}(M, F)\neq\emptyset$. So, we immediately have the following proposition.

\begin{pn}\label{pn-sup-str-sL} Ultra sL-approximation spaces are all
sL-approximation spaces.
\end{pn}

Next, we explore properties of
sL-approximation spaces.
\begin{pn}\label{3-pn-M-sup}
Let $(U, R, \mathcal{F})$ be an sL-approximation space, $ F_1, F_2\in \mathcal{F}$, $M, M_{1}, M_{2}\in \mathcal{F}^{\sharp}$ and $E\in \mathfrak{C}(U, R, \mathcal{F})$.  Then the following statements hold.

$(1)$ If    $F_{1}\cup F_{2}\subseteq E$, $M\subseteq\overline{R}(F_{1})\cap\overline{R}(F_{2})$, then $\mathcal{S}(M, F_{1})=\mathcal{S}(M, F_{2})$.

$(2)$ If    $F_{1}\cup F_{2}\subseteq E$,  $M_{1}\subseteq\overline{R}(F_1)$, $M_{2}\subseteq\overline{R}(F_{2})$,  $\overline{R}(M_{1})\subseteq\overline{R}(M_{2})$
and $G_{i}\in
\mathcal{S}(M_{i}, F_i)\ (i=1, 2)$, then $\overline{R}(G_{1})\subseteq\overline{R}(G_{2})$.
\end{pn}
\noindent{\bf Proof.}
(1) By $F_{1}\cup F_{2}\subseteq_{fin} E\in \mathfrak{C}(U, R, \mathcal{F})$ and Definition \ref{dn-cf-cl},   there exists $F_{3} \in \mathcal{F}$ such that $F_{1}\cup F_{2}\subseteq \overline{R}(F_{3})\subseteq E$. It follows from Lemma \ref{lm-ref-lm2-tr-up}(3) that $\overline{R}(F_{1})\cup\overline{R}(F_{2})\subseteq\overline{R}(F_{3})$.
Thus $M\subseteq\overline{R}(F_{1})\subseteq\overline{R}(F_{3})$ and $M\subseteq\overline{R}(F_{2})\subseteq\overline{R}(F_{3})$. So, $\mathcal{S}(M, F_{i})\ne\emptyset\ (i=1,2,3)$.
 By Proposition \ref{1-pn-M-sup}(2), we have $\mathcal{S}(M, F_{1})=\mathcal{S}(M, F_{3})$ and $\mathcal{S}(M, F_{2})=\mathcal{S}(M, F_{3})$. Thus, $\mathcal{S}(M, F_{1})=\mathcal{S}(M, F_{2})$.

(2) Clearly,  there exists $F_{3} \in \mathcal{F}$ such that $F_{1}\cup F_{2}\subseteq \overline{R}(F_{3})\subseteq E$. So, we have $M_1\subseteq\overline{R}(F_{1})\subseteq\overline{R}(F_{3})$ and $M_2\subseteq\overline{R}(F_{2})\subseteq\overline{R}(F_{3})$ and $\mathcal{S}(M_i, F_{i})\ne\emptyset\ne\mathcal{S}(M_i, F_{3})\ (i=1,2)$.  By Proposition \ref{1-pn-M-sup}(2), we have $\mathcal{S}(M_{1}, F_{1}) =\mathcal{S}(M_{1}, F_{3})$ and $\mathcal{S}(M_{2}, F_{2}) =\mathcal{S}(M_{2}, F_{3})$. It follows from $G_{i}\in
\mathcal{S}(M_{i}, F_i)\ (i=1, 2)$ that $G_{1}\in \mathcal{S}(M_{1}, F_{3})$ and $G_{2}\in \mathcal{S}(M_{2}, F_{3})$. By   $\overline{R}(M_{1})\subseteq\overline{R}(M_{2})$ and Proposition \ref{1-pn-M-sup}(3),  we have $\overline{R}(G_{1})\subseteq\overline{R}(G_{2})$.\hfill$\Box$

\begin{tm}\label{tm-sL-sLDO} If $(U, R, \mathcal{F})$ is an sL-approximation space, then the dcpo $(\mathfrak{C}(U, R, \mathcal{F}), \subseteq)$ is an sL-domain.
\end{tm}
\noindent{\bf Proof.}   Clearly,  $(\mathfrak{C}(U, R, \mathcal{F}), \subseteq)$ is a continuous domain. So it suffices to show that for every$E\in (\mathfrak{C}(U, R, \mathcal{F}), \subseteq)$, the set ${\downarrow\!E}=\{H\mid H\in \mathfrak{C}(U, R, \mathcal{F}), H\subseteq E\}$ is a sup-semilattice ordered by set inclusion. For any $H_{1}, H_{2}\in {\downarrow\!E}$, assume that $H_{1}\cup H_{2}\neq\emptyset$ without losing generality.  Construct
\vskip 3pt
\centerline{$\mathcal{H}=\bigcup\{\mathcal{S}(F_{1}\cup F_{2}, F_{3})\mid  F_1,F_2,F_3\in \mathcal{F}, F_{i}\subseteq H_{i}\ (i=1,2), F_{3}\subseteq E, F_1\cup F_2\subseteq \overline{R}(F_3)  \}$}
\vskip 3pt
\noindent and
\vskip 3pt
\centerline{$H=\bigcup\{\overline{R}(G)\mid G\in \mathcal{H}\}$.}
\vskip 3pt
\noindent To show $H$ is the least upper bound of $H_{1}$ and $H_{2}$ in ${\downarrow\!E}$, we divide the  proof into several steps.

Step 1. It follows directly from Lemma \ref{pn2-cf-clo}(2)  that $ H\subseteq  E$.

Step 2. Assert that $ H_{1}\cup H_{2}\subseteq H$.

Let $x\in  H_{1}\cup H_{2}$. Assume that $x\in  H_{1}$ without losing generality. Then by Definition \ref{dn-cf-cl}, there exists $F_{1}\in \mathcal{F}$ such that $x\in \overline{R}(F_{1})$ and $F_1\subseteq H_{1}$. Select $F_{2}\in \mathcal{F}$ with $F_{2}\subseteq  H_{2}$. Then $F_1\cup F_2\subseteq H_{1}\cup H_{2}\subseteq E$. It follows from $E\in \mathfrak{C}(U, R, \mathcal{F})$ that there exists $F_{3}\in \mathcal{F}$ such that $F_1\cup F_2\subseteq \overline{R}(F_3)$ and $F_3\subseteq E$.  By the definition of sL-approximation spaces, we have $\mathcal{S}(F_1\cup F_2, F_3)\ne\emptyset$. For $G\in \mathcal{S}(F_1\cup F_2, F_3)\subseteq\mathcal{H}$, we have $x\in \overline{R}(F_1\cup F_2)=\overline{R}(F_1)\cup \overline{R}(F_2)\subseteq\overline{R}(G)\subseteq H$. By arbitrariness of $x\in  H_{1}\cup H_{2}$, we have  $ H_{1}\cup H_{2}\subseteq H$.

Step 3. Assert that $ H\in \mathfrak{C}(U, R, \mathcal{F})$.

It suffices to show  the family $\{\overline{R}(G)\mid G\in \mathcal{H}\}$ is directed. Clearly, $\{\overline{R}(G)\mid G\in \mathcal{H}\}\neq\emptyset$. Let $G_{1}\in \mathcal{H}$. Then there exist $F_{11}, F_{12}, F_{13}\in \mathcal{F}$  such that
\vskip3pt

\centerline{$ F_{11}\subseteq H_{1},   F_{12}\subseteq H_2, F_{13}\subseteq E, F_{11}\cup F_{12}\subseteq \overline{R}(F_{13}) \mbox{ and } G_1\in \mathcal{S}(F_{11}\cup F_{12}, F_{13}).$}
\vskip 5pt

\noindent Let $ G_{2}\in \mathcal{H}$. Then there exist $F_{21}, F_{22}, F_{23}\in \mathcal{F}$ such that
\vskip3pt

\centerline{$ F_{21}\subseteq H_{1},   F_{22}\subseteq H_2, F_{23}\subseteq E, F_{21}\cup F_{22}\subseteq \overline{R}(F_{23})\mbox{ and }G_2\in \mathcal{S}(F_{21}\cup F_{22}, F_{23}).$}
\vskip5pt
\noindent  It follows from $F_{11}\cup F_{21}\subseteq_{fin} H_1\in  \mathfrak{C}(U, R, \mathcal{F})$ that there exists $F_1\in \mathcal{F}$ such that $F_{11}\cup F_{21}\subseteq \overline{R}(F_1)$ and $F_{1}\subseteq H_1$. Similarly, there exists $F_2\in \mathcal{F}$ such that $F_{12}\cup F_{22}\subseteq \overline{R}(F_2)$ and $F_{2}\subseteq H_2$. Thus, $F_{11}\cup F_{12}\subseteq \overline{R}(F_1\cup F_2)$ and $F_{21}\cup F_{22}\subseteq \overline{R}(F_1\cup F_2)$. Since $F_{1}\cup F_{2}\subseteq H_1\cup H_2 \subseteq E$, there exists $F_3\in \mathcal{F}$ such that $F_{1}\cup F_{2}\subseteq \overline{R}(F_3)$ and $F_3\subseteq E$. Let $G_3\in \mathcal{S}(F_{1}\cup F_{2}, F_{3})\subseteq \mathcal{H}$. Noticing that $F_{3}\cup F_{13}\subseteq_{fin} E$ and $F_{11}\cup F_{12}\subseteq \overline{R}(F_1\cup F_2)$, by Proposition \ref{3-pn-M-sup}(2), we have $\overline{R}(G_{1})\subseteq\overline{R}(G_{3})$. Similarly,
 $\overline{R}(G_{2})\subseteq\overline{R}(G_{3})$. Thus  $(\{\overline{R}(G)\mid G\in \mathcal{H}\}, \subseteq)$ is directed. It follows from Lemma \ref{pn3-cf-cl}(3) that $H=\bigcup\{\overline{R}(G)\mid G\in \mathcal{H}\}\in\mathfrak{C}(U, R, \mathcal{F})$.

Step 4. Assert that $H$ is the least upper bound of $H_{1}$ and $H_{2}$ in ${\downarrow\!E}$.

Let $H^{\prime}\in {\downarrow\!E}$ and $H_{1}\cup H_{2}\subseteq H^{\prime}$. For all $x\in H$, there exist $F_{i}\in\mathcal{F}\ (i=1,2,3)$ such that
\vskip3pt
\centerline{$F_{i}\subseteq H_{i}\  (i=1,2), F_{3}\subseteq E, F_{1}\cup F_{2}\subseteq \overline{R}(F_{3}), G\in \mathcal{S}(F_{1}\cup F_{2}, F_{3})$, $x\in \overline{R}(G)$.}
\vskip5pt
\noindent It follows from $F_{1}\cup F_{2}\subseteq_{fin}
H_{1}\cup H_{2}\subseteq H^{\prime}$ and $H^{\prime}\in \mathfrak{C}(U, R, \mathcal{F})$ that there exists $F_{4}\in \mathcal{F}$ such that $F_{1}\cup F_{2}\subseteq \overline{R}(F_{4})\subseteq H^{\prime}$ and $F_{4}\subseteq H^{\prime}$.   Noticing that $F_3\cup F_4\subseteq E$,  by proposition \ref{3-pn-M-sup}(1), we have $G\in \mathcal{S}(F_{1}\cup F_{2}, F_{3})=\mathcal{S}(F_{1}\cup F_{2}, F_{4})$. Thus $x\in \overline{R}(G)\subseteq\overline{R}(F_{4})\subseteq H^{\prime}$. By  arbitrariness of $x\in  H$, we have $H\subseteq H^{\prime}$, showing that $H$ is the least upper bound of $H_{1}$ and $H_{2}$ in ${\downarrow\!E}$.

To sum up, $(\mathfrak{C}(U, R, \mathcal{F}), \subseteq)$ is an sL-domain.\hfill$\Box$\\

The reader may wonder whether   dense  subspaces of  sL-approximation spaces are sL-approximation spaces.  For topological CF-approximation spaces, we have the following affirmative answer.

 \begin{pn} Let  $(U, R, \mathcal{F})$ be  a topological  sL-approximation space, $(V, R|_V, \mathcal{G})$ a dense CF-approximation subspace of $(U, R, \mathcal{F})$. Then $(V, R|_V, \mathcal{G})$ is an sL-approximation space.
\end{pn}
\noindent{\bf Proof.}
For all $M\in \mathcal{G}^{\sharp}$ and $F\in \mathcal{G}$ with $M\subseteq \overline{R|_V}(F)$, we have that $M\subseteq \overline{R}(F)$  and $\overline{R}(M)\subseteq \overline{R}(F)$ hold in $(U, R, \mathcal{F})$.
Since $(U, R, \mathcal{F})$ is an sL-approximation space, we have that $\mathcal{S}(M, F)\ne \emptyset$. Let $G_1\in \mathcal{S}(M, F)$. Then by
(sL-1) and $R$ being a preorder, we have $M\subseteq\overline{R}(M)\subseteq \overline{R}(G_1)\subseteq \overline{R}(F)$. By the density of $(V, R|_V, \mathcal{G})$, there exists   a $G\in \mathcal{G}$ such that  $M \subseteq \overline{R}(G)\subseteq \overline{R}(G_1)\subseteq\overline{R}(F)$.  Since $G_1$ satisfies (sL-2), we have that $\overline{R}(G_1)\subseteq \overline{R}(G)$ and thus $\overline{R}(G_1)=\overline{R}(G)$. By Proposition \ref{1-pn-M-sup}(1) that $G\in \mathcal{S}(M, F)$.

 Next we show  that $G$ satisfies (sL-1) and  (sL-2) in $(V, R|_V, \mathcal{G})$. Clearly, (sL-1) $\overline{R|_V}(M)\subseteq \overline{R|_V}(G)\subseteq \overline{R|_V}(F)$  holds in $(V, R|_V, \mathcal{G})$. To show (sL-2), we need show for all $G^{\prime}\in \mathcal{G}$, $\overline{R|_V}(M)\subseteq \overline{R|_V}(G^{\prime})\subseteq \overline{R|_V}(F)$ implies that $\overline{R|_V}(G)\subseteq\overline{R|_V}(G^{\prime})$ also holds in $(V, R|_V, \mathcal{G})$.
It follows from  $\overline{R|_V}(M)\subseteq \overline{R|_V}(G^{\prime})=\overline{R}(G^{\prime})\cap V$ and $R$ being a preorder that  $M\subseteq \overline{R}(G^{\prime})$ and $\overline{R}(M)\subseteq \overline{R}(G^{\prime})$. It follows from $G^{\prime}\subseteq \overline{R|_V}(G^{\prime})\subseteq \overline{R|_V}(F)=\overline{R}(F)\cap V$  that  $G^{\prime}\subseteq \overline{R}(F)$. Thus we have that $\overline{R}(M)\subseteq \overline{R}(G^{\prime})\subseteq \overline{R}(F)$. By $G$ satisfying (sL-2) in $(U, R, \mathcal{F})$, we have $\overline{R}(G)\subseteq\overline{R}(G^{\prime})$ and $\overline{R|_V}(G)\subseteq\overline{R|_V}(G^{\prime})$, showing that (sL-2) hold in $(V, R|_V, \mathcal{G})$.  So, by Definition \ref{dn-sL-App},   $(V, R|_V, \mathcal{G})$ is an sL-approximation space.\hfill$\Box$\\


Next, we consider a special  kind of sL-approximation spaces---L-approximation spaces.
\begin{dn}
Let $(U, R, \mathcal{F})$ be a CF-approximation space.  If $\mathcal{S}(M, F)\neq\emptyset$ for all
 $M\in \mathcal{F}^{\sharp}\cup\{\emptyset\}$ with $M\subseteq\overline{R}(F)$;  that is to say, there exists $G\in \mathcal{F}$ satisfies the following conditions for all
 $M\in \mathcal{F}^{\sharp}\cup\{\emptyset\}$ with $M\subseteq\overline{R}(F)$:

{(sL-1)} $\overline{R}(M)\subseteq\overline{R}(G)\subseteq\overline{R}(F)$; and

{ (sL-2)} $\forall G^{\prime} \in \mathcal{F}$, $\overline{R}(M)\subseteq\overline{R}(G^{\prime})\subseteq\overline{R}(F)\Longrightarrow\overline{R}(G)\subseteq\overline{R}(G^{\prime})$,

\noindent then $(U, R, \mathcal{F})$ is called an {\em L-approximation space}.
\end{dn}
\begin{rk}

 Clearly,  by  shifting the precondition ``$M\in \mathcal{F}^{\sharp}$" in  Definition \ref{dn-sL-App} of sL-approximation spaces to ``$M\in
\mathcal{F}^{\sharp}\cup\{\emptyset\}$",  we get the concept of L-approximation spaces. Thus  L-approximation spaces are all sL-approximation spaces.
For an sL-approximation space $(U, R, \mathcal{F})$, if $\emptyset\in \mathcal{F}$, then $(U, R, \mathcal{F})$ is an L-approximation space automatically.

\end{rk}

%
%
%
%
%
%
\begin{tm}\label{tm-L-LDO} For an L-approximation space $(U, R, \mathcal{F})$, the dcpo $(\mathfrak{C}(U, R, \mathcal{F}), \subseteq)$ is an L-domain.
\end{tm}
\noindent{\bf Proof.}   It follows from Theorem \ref{tm-sL-sLDO} that  $(\mathfrak{C}(U, R, \mathcal{F}), \subseteq)$ is an sL-domain. To show  $(\mathfrak{C}(U, R, \mathcal{F}), \subseteq)$ is an L-domain, it suffices to prove for all  $E\in (\mathfrak{C}(U, R, \mathcal{F}), \subseteq)$, the principal ideal $\downarrow\! E$  has a least element. It follows from $E\in \mathfrak{C}(U, R, \mathcal{F})$ and $\emptyset\subseteq_{fin}E$ that there exists $F\in \mathcal{F}$ such that $F\subseteq E$ and $\emptyset\subseteq \overline{R}(F)$. Let $G\in \mathcal{S}(\emptyset, F)$,  we next show $\overline{R}(G)$ is precisely the least element of $(\downarrow\! E, \subseteq)$. Let $H\in \downarrow\! E$. Then there exists $F_1\in \mathcal{F}$ such that $F_1\subseteq H\subseteq E$. It follows from $F\cup F_1\subseteq_{fin} E$ and Proposition \ref{3-pn-M-sup}(1) that  $ \mathcal{S}(\emptyset, F)=\mathcal{S}(\emptyset, F_1)$. Thus, $G\in \mathcal{S}(\emptyset, F_1)$ and $\overline{R}(G)\subseteq\overline{R}(F_1)\subseteq H$, showing that $\overline{R}(G)$ is the  least element of $\downarrow\! E$. By   arbitrariness of $E\in \mathfrak{C}(U, R, \mathcal{F})$, we have that $(\mathfrak{C}(U, R, \mathcal{F}), \subseteq)$ is an L-domain.\hfill$\Box$

\begin{tm} \label{rep-sL-topsL} For a topological sL-approximation space (resp.,  topological L-approximation space) $(U, R, \mathcal{F})$, the dcpo  $(\mathfrak{C}(U, R, \mathcal{F}), \subseteq)$ is an algebraic  sL-domain (resp., algegbraic L-domain).
\end{tm}
\noindent{\bf Proof.}  It follows from Lemma \ref{lm-top-alg} that $(\mathfrak{C}(U, R, \mathcal{F}), \subseteq)$ is an algebraic domain. By Theorem \ref{tm-sL-sLDO} (resp., Theorem \ref{tm-L-LDO}), we know that $(\mathfrak{C}(U, R, \mathcal{F}), \subseteq)$ is an sL-domain (resp.,  L-domain), and hence an algebraic sL-domain (resp.,  L-domain).
 \hfill$\Box$
\begin{rk}\label{rk-top-equ} Let $(U, R, \mathcal{F})$ be  a topological CF-approximation space. For all $M\in \mathcal{F}^{\sharp}, F\in \mathcal{F}$,   by Lemma \ref{lm-ref-lm2-tr-up}(4), we have that $M\subseteq \overline{R}(F)\Leftrightarrow\overline{R}(M)\subseteq \overline{R}(F)$, (sL-1)$\Leftrightarrow$ (sL-1$^{\prime})$ and $\mathcal{S}(M, F)=\mathcal{S}^{*}(M, F)$. So, $(U, R, \mathcal{F})$ is a topological ultra sL-approximation space if and only if $(U, R, \mathcal{F})$
is a  topological sL-approximation space.
\end{rk}

\section{Representations of  sL-domains}
In this section, by means with sL-approximation spaces (resp., L-approximation spaces, bc-approximation spaces), we  shall discuss  representations for sL-domains (resp., L-domains, bc-domains).

\begin{tm}\label{tm-CDO-sLapp} Let  $(L, R_{L},  \mathcal{F}_{L})$ be  the induced CF-approximation space by  an sL-domain  $(L, \leqslant)$. Then
 $(L, R_{L},  \mathcal{F}_{L})$ is  an sL-approximation space.
\end{tm}
\noindent{\bf Proof.}
To prove that $(L, R_{L},  \mathcal{F}_{L})$ is an  sL-approximation space, let $M\in (\mathcal{F}_{L})^{\sharp}$, $F\in \mathcal{F}_{L}$ and $M\subseteq \overline{R_L}(F)$. We use $c_F$ to denote the greatest element in $F$. Thus $M \subseteq_{fin}\overline{R_L}(F)=\dda c_F$. Since $L$ is an sL-domain, we have that $\bigvee_{c_F}M$ exists.  Set $G=\{\bigvee_{c_F}M\}\in \mathcal{F}_{L}$. Clearly, $\overline{R_{L}}(G)=\dda \bigvee_{c_F}M\subseteq\dda c_F=\overline{R_{L}}(F)$  and $G$ satisfies (sL-1).
To check that $G$ satisfies (sL-2), let $H\in \mathcal{F}_L$ and $\overline{R_{L}}(M)\subseteq\overline{R_{L}}(H)\subseteq\overline{R_{L}}(F)$. It follows from Lemma \ref{lm-ula} that for all $m\in M$, $\dda m\subseteq \dda c_{H}\subseteq \dda c_{F}$. Since $L$ is a continuous domain, we have $M\subseteq\downarrow\! c_{H}\subseteq\downarrow\! c_{F}$. By Lemma \ref{lm-sup-st}, we have $\bigvee_{c_H}M=\bigvee_{c_F}M\leqslant c_H$. Thus $\overline{R_{L}}(G)=\dda \bigvee_{c_F}M\subseteq\dda  c_H=\overline{R_{L}}(H)$, showing that $G$ satisfies (sL-2). Thus $G\in \mathcal{S}(M, F)\neq\emptyset$. This shows that $(L, R_{L},  \mathcal{F}_{L})$ is an  sL-approximation space.
\hfill$\Box$

\begin{ex}\label{ex-sl<ultra}
Let $L=[0, 10]$ be a closed interval with the ordinary order of number. Clearly, $L$ is an sL-domain. By Theorem \ref{tm-CDO-sLapp}, the  induced CF-approximation space $(L, R_L, \mathcal{F}_{L})$ is an  sL-approximation space. However, $(L, R_L, \mathcal{F}_{L})$ is not an ultra sL-approximation space. To see this, assume $(L, R_L, \mathcal{F}_{L})$ is an ultra sL-approximation space. Then for $M=\{2, 4\}\in \mathcal{F}_{L}^{\sharp}$ and $F=\{1, 3, 4\}\in \mathcal{F}_{L}$, we have $\overline{R_{L}}(M)=[0, 4)\subseteq \overline{R_{L}}(F)=[0, 4)$. By the assumption, there exists $G\in \mathcal{F}_{L}$ satisfies (sL-1$^{\prime}$), that is $\overline{R_{L}}(M)=[0, 4)\subseteq\overline{R_{L}}(G)=[0, c_G)$ and $G\subseteq\overline{R_{L}}(F)=[0, 4)$. It would be that  $c_G\geqslant 4$ and $c_{G}< 4$, a contradiction. This shows that  the sL-approximation space $(L, R_L, \mathcal{F}_{L})$ is not an ultra sL-approximation space.
\end{ex}
%
%


\begin{tm}\label{tm-rep-sL} {\em (Representation Theorem  for sL-domains)} For a poset $(L, \leqslant)$,  the following statements are equivalent.

 $(1)$ $L$ is an  sL-domain;

 $(2)$ There exists an sL-approximation space $(U, R, \mathcal{F})$ such that $(L, \leqslant)\cong (\mathfrak{C}(U, R, \mathcal{F}), \subseteq)$.

%
\end{tm}
\noindent{\bf Proof.}
$(1)\Rightarrow (2)$: If $L$ is an sL-domain, then by Theorem \ref{tm-CDO-sLapp},  the induced CF-approximation space $(L, R_{L},  \mathcal{F}_{L})$ is an sL-approximation space. Define  $f: L\to \mathfrak{C}(L, R_{L},  \mathcal{F}_{L})$ such that for all $x\in L$, $f(x)=\dda x$. Then it follows from  Example \ref{tm-CDO-CFGA}(2) and the continuity of $L$ that $f$ is an order-isomorphism.


$(2)\Rightarrow (1)$:   Follows directly from Theorem \ref{tm-sL-sLDO}.
  \hfill$\Box$

\begin{tm}\label{tm-rep-Ldo} {\em (Representation Theorem  for L-domains)} Let $(L, \leqslant)$ be a poset. Then $L$ is an  L-domain iff there exists an L-approximation space  $(U, R, \mathcal{F})$ such that $(L, \leqslant)\cong (\mathfrak{C}(U, R, \mathcal{F}), \subseteq)$.
\end{tm}
\noindent{\bf Proof.}
$\Leftarrow$: Follows directly from Theorem \ref{tm-L-LDO}.

$\Rightarrow$: Let $L$ be an  L-domain. It suffices to prove that $(L, R_{L},  \mathcal{F}_{L})$ is an L-approximation space. It follows from Theorem \ref{tm-CDO-sLapp} that $(L, R_{L},  \mathcal{F}_{L})$ is an  sL-approximation space . To prove that $(L, R_{L},  \mathcal{F}_{L})$ is an L-approximation space, let $M=\emptyset$ and $F\in \mathcal{F}^{\sharp}$. Since $L$ is an L-domain, we have that ${\downarrow} c_{F}$ has the least element denoted by $\bot_{F}$. Set $G=\{ \bot_{F}\}\in \mathcal{F}_{L}$. Clearly, $\overline{R_{L}}(\emptyset)=\emptyset\subseteq\overline{R_{L}}(G)=\dda \bot_{F}\subseteq\dda {c_F}=\overline{R_{L}}(F)$, showing that $G$ satisfies (sL-1). If $G_{1}\in \mathcal{F}$ and $\overline{R_{L}}(\emptyset)=\emptyset\subseteq\overline{R_{L}}(G_{1})=\dda c_{G_{1}}\subseteq\dda c_F=\overline{R_{L}}(F)$, then by the continuity of $L$, we have $c_{G_1}\leqslant c_F$. Since  $\bot_{F}$ is the least element in ${\downarrow} c_{F}$, we have $\bot_{F}\leqslant c_{G_{1}}$. Therefore $\overline{R_{L}}(G)=\dda \bot_{F}\subseteq\dda c_{G_{1}}=\overline{R_{L}}(G_{1})$, showing that $G$ satisfies (sL-2). Thus  $(L, R_{L},  \mathcal{F}_{L})$ is an  L-approximation space.
  \hfill$\Box$

\begin{dn}
Let $(U, R, \mathcal{F})$ be a CF-approximation space.  If for all $M\in \mathcal{F}^{\sharp}\cup\{\emptyset\}$ and $F\in \mathcal{F}$ with $M\subseteq\overline{R}(F)$, there exists $G\in \mathcal{F}$ satisfies the following conditions:

{(sL-$1)$} $\overline{R}(M)\subseteq\overline{R}(G)\subseteq\overline{R}(F)$;

{ (sL-$2^{\prime})$} $\forall G^{\prime} \in \mathcal{F}$, $\overline{R}(M)\subseteq\overline{R}(G^{\prime})\Longrightarrow\overline{R}(G)\subseteq\overline{R}(G^{\prime})$,

\noindent then $(U, R, \mathcal{F})$ is called a {\em bc-approximation space}.
\end{dn}
\begin{rk} Clearly, (sL-$2^{\prime})$ implies (sL-2). So bc-approximation spaces are all L-approximation spaces. For all $M\in \mathcal{F}^{\sharp}\cup\{\emptyset\}$ and $F\in \mathcal{F}$ with $M\subseteq\overline{R}(F)$, we use $\Sigma(M, F)$ to denote all $G\in \mathcal{F}$ that satisfies (sL-1) and (sL-$2^{\prime})$. Obviously, $\Sigma(M, F)\subseteq \mathcal{S}(M, F)$.
\end{rk}
\begin{pn}\label{pn-bc-sup}
Let $(U, R, \mathcal{F})$ be a bc-approximation space, $F_1, F_2\in \mathcal{F}$ and $M, M_{1}, M_{2}\in \mathcal{F}^{\sharp}\cup\{\emptyset\}$.  Then the following statements hold.

$(1)$ If     $M\subseteq\overline{R}(F)$, $G_{1}\in \Sigma(M, F)$ and $G_{2}\in \mathcal{F}$, then $G_{2}\in \Sigma(M, F)$ iff $\overline{R}(G_{1})=\overline{R}(G_{2})$.

$(2)$ If    $M\subseteq\overline{R}(F_{1})\cap\overline{R}(F_{2})$, then $\Sigma(M, F_{1})=\Sigma(M, F_{2})$.

$(3)$ If     $M_{1}\subseteq\overline{R}(F_1)$,
$\overline{R}(M_{1})\subseteq\overline{R}(M_{2})$, $M_{2}\subseteq\overline{R}(F_{2})$,  $G_{1}\in \Sigma(M_{1}, F_{1})$ and $G_{2}\in \Sigma(M_{2}, F_{2})$, then $\overline{R}(G_{1})\subseteq\overline{R}(G_{2})$.
\end{pn}
\noindent{\bf Proof.}
(1) Similar to the proof of Proposition \ref{1-pn-M-sup}(1).

(2) Let $G_{1}\in \Sigma(M, F_{1})$ and $G_{2}\in \Sigma(M, F_{2})$. By (sL-1), we have $\overline{R}(M)\subseteq\overline{R}(G_{1})\subseteq\overline{R}(F_{1})$ and $\overline{R}(M)\subseteq\overline{R}(G_{2})\subseteq\overline{R}(F_{2})$. By (sL-$2^{\prime})$, we have $\overline{R}(G_{1})\subseteq\overline{R}(G_{2})$ and $\overline{R}(G_{2})\subseteq\overline{R}(G_{1})$. Thus $\overline{R}(G_{1})=\overline{R}(G_{2})$. By (1),  we have $\Sigma(M, F_{1})=\Sigma(M, F_{2})$.

(3)
It follows from $G_{2}\in \Sigma(M_{2}, F_{2})$ and $\overline{R}(M_{1})\subseteq\overline{R}(M_{2})$ that $\overline{R}(M_{1})\subseteq\overline{R}(M_{2})\subseteq\overline{R}(G_{2})$. By  $G_{1}\in \Sigma(M_{1}, F_{1})$ and (sL-2$^{\prime}$) that $\overline{R}(G_{1})\subseteq\overline{R}(G_{2})$.
\hfill$\Box$
\begin{tm}\label{tm-bc-bcdom}
Let $(U, R, \mathcal{F})$ be a bc-approximation space.  Then $(\mathfrak{C}(U, R, \mathcal{F}), \subseteq)$ is a bc-domain.
\end{tm}
\noindent{\bf Proof.}  Clearly, $(\mathfrak{C}(U, R, \mathcal{F}),\subseteq)$ is a continuous domain. So it suffices to show  $(\mathfrak{C}(U, R, \mathcal{F}), \subseteq)$ has the least element and any up-bounded two-element family $\{H_{1}, H_{2}\}\subseteq (\mathfrak{C}(U, R, \mathcal{F}), \subseteq)$ has  the supremum in $(\mathfrak{C}(U, R, \mathcal{F}), \subseteq)$.

Let $F\in \mathcal{F}$ and $G\in \Sigma(\emptyset, F)\ne\emptyset$. It  follows from Proposition \ref{pn-bc-sup}(2) that $\overline{R}(G)$  is the least element in $(\mathfrak{C}(U, R, \mathcal{F}), \subseteq)$.

Let $E$ be an upper bound of $\{H_{1}, H_{2}\}$. Construct
\vskip 3pt
\centerline{$\mathcal{H}=\bigcup\{ \Sigma(F_{1}\cup F_{2}, F_{3})\in \mathcal{F}\mid  F_{i}\in \mathcal{F}\ (i=1,2,3), F_{i}\subseteq H_{i}\ (i=1,2), F_{3}\subseteq E, F_1\cup F_2\subseteq \overline{R}(F_3)\}$}
\vskip 3pt
\noindent and
\vskip 3pt
\centerline{$H=\bigcup\{\overline{R}(G)\mid G\in \mathcal{H}\}$.}
\vskip 3pt
\noindent We divide the  proof of $H$ being the supremum of $\{H_{1}, H_{2}\}$ in $(\mathfrak{C}(U, R, \mathcal{F}),\subseteq)$ into several steps.

Step 1. It follows directly from Lemma \ref{pn2-cf-clo}(2) that $H\subseteq  E$.

 Step 2.
Similar to the proof of Step 2  in Theorem \ref{tm-sL-sLDO}, we  have  $ H_{1}\cup H_{2}\subseteq H$.

Step 3. Similar to the proof of Step 3 in Theorem  \ref{tm-sL-sLDO}, we have $ H\in \mathfrak{C}(U, R, \mathcal{F})$.

Step 4. Assert that $H$ is the least upper bound of $H_{1}$ and $H_{2}$ in $(\mathfrak{C}(U, R, \mathcal{F}), \subseteq)$.

Let $H^{\prime}\in \mathfrak{C}(U, R, \mathcal{F})$ and $H_{1}\cup H_{2}\subseteq H^{\prime}$. For all $x\in H$, there exist $F_{i}\in\mathcal{F}\ (i=1,2,3)$ such that
\vskip3pt
\centerline{$F_{i}\subseteq H_{i}\ (i=1,2), F_{3}\subseteq E$, $F_1\cup F_2\subseteq \overline{R}(F_3)$, $G\in \Sigma(F_{1}\cup F_{2}, F_{3})$ and $x\in \overline{R}(G)$.}
\vskip5pt
\noindent  It follows from $F_{1}\cup F_{2}\subseteq_{fin}
H_{1}\cup H_{2}\subseteq H^{\prime}\in \mathfrak{C}(U, R, \mathcal{F})$ that there exists $F_{4}\in \mathcal{F}$ such that $F_{1}\cup F_{2}\subseteq \overline{R}(F_{4})\subseteq H^{\prime}$ and $F_{4}\subseteq H^{\prime}$. By Proposition \ref{pn-bc-sup}(2), we have $G\in \Sigma(F_{1}\cup F_{2}, F_{3})=\Sigma(F_{1}\cup F_{2}, F_{4})$. Thus $x\in \overline{R}(G)\subseteq\overline{R}(F_{4})\subseteq H^{\prime}$. By  arbitrariness of $x\in  H$, we have $H\subseteq H^{\prime}$, showing that $H$ is the least upper bound of $H_{1}$ and $H_{2}$ in $(\mathfrak{C}(U, R, \mathcal{F}), \subseteq)$.

To sum up,  we see that $(\mathfrak{C}(U, R, \mathcal{F}), \subseteq)$ is a bc-domain.\hfill$\Box$
\begin{tm}\label{th-rep-bc-domain} {\em (Representation Theorem for bc-domains)} A poset $(L, \leqslant)$ is a bc-domain if and only if there exists some bc-approximation space such that $(\mathfrak{C}(U, R, \mathcal{F}), \subseteq)\cong (L, \leqslant)$.
\end{tm}
\noindent\noindent{\bf Proof.}
$\Leftarrow$: Follows from Theorem \ref{tm-bc-bcdom}.

$\Rightarrow$: It suffices to prove that for  a bc-domain $L$, $(L, R_{L},  \mathcal{F}_{L})$ is a bc-approximation space. Let $M\in (\mathcal{F}_{L})^{\sharp}\cup\{\emptyset\}$, $F\in \mathcal{F}_{L}$ with $M\subseteq \overline{R_L}(F)$.  We use $c_F$ to denote the greatest element in $F$. If $M =\emptyset$, let $\bot$ be the least element of $L$ and $G=\{\bot\}\in\mathcal{F}_{L}$,  then $G\in \Sigma(\emptyset, F)$. If $M\ne\emptyset$,
then $M\subseteq_{fin}\overline{R_L}(F)=\dda c_F$. Since $L$ is a bc-domain, we have that $\bigvee M$ exists.
Set $G=\{\bigvee M\}\in \mathcal{F}_L$. Clearly,  $G\in \Sigma(M, F)$, that is $G$ satisfies (sL-1) and  {(sL-$2^{\prime})$}.  Thus  $(L, R_{L},  \mathcal{F}_{L})$ is a bc-approximation space.\hfill$\Box$
\begin{pn}\label{pn2-hy-sL} A topological CF-approximation space $(U, R, \mathcal{F})$ is an sL-approximation space if and only if $(\{\overline{R}(F)\mid F\in \mathcal{F}\}, \subseteq)$ is an  sL-cusl.
\end{pn}
\noindent {\bf Proof.} 
$\Longrightarrow$:  Let $F, F_{1}, F_{2}\in \mathcal{F}$ and $\overline{R}(F_{1})\cup \overline{R}(F_{2})\subseteq\overline{R}(F)$. By Lemma \ref{lm-ula} and \ref{lm-ref-lm2-tr-up}(4),
we have
$F_{1}\cup F_{2}\subseteq\overline{R}(F_{1}\cup F_{2})\subseteq\overline{R}(F)$. It follows from $F_{1}\cup F_{2}\in F^{\sharp}$ that there exists $G\in \mathcal{S}(F_{1}\cup F_{1}, F)$ such that $G$ satisfies (sL-1) and (sL-2). By (sL-1), we have that $\overline{R}(G)$ is an upper bound of $ \overline{R}(F_{1}),  \overline{R}(F_{2})$ in

\vskip 3pt
\centerline{${\downarrow}\overline{R}(F)=\{\overline{R}(G)\mid G\in \mathcal{F}, \overline{R}(G)\subseteq\overline{R}(F)\}$.}
\vskip5pt

\noindent By (sL-2), we have that  $\overline{R}(G)$ is the least upper bound of $ \overline{R}(F_{1}),  \overline{R}(F_{2})$ in $\downarrow\!\overline{R}(F)$, showing that $\downarrow\!\overline{R}(F)$ is a sup-semilattice. Thus $(\{\overline{R}(F)\mid F\in \mathcal{F}\}, \subseteq)$ is an  sL-cusl.

$\Longleftarrow$: Let $M=\bigcup_{i=1}^{n}F_{i}\in \mathcal{F}^{\sharp}$ and $F\in \mathcal{F}$  with $M\subseteq\overline{R}(F)$. By Lemma \ref{lm-ula}, we have
$\overline{R}(M)=\bigcup_{i=1}^{n}\overline{R}(F_{i})$. Thus $\{\overline{R}(F_{i})\mid i\leqslant n\}\subseteq\downarrow\! \overline{R}(F)$. It follows from $(\{\overline{R}(G)\mid G\in \mathcal{F}\}, \subseteq)$ is an  sL-cusl that there exists $G\in \mathcal{F}$ such that $\overline{R}(G)$ is the least upper bound of $\{\overline{R}(F_{i})\mid i\leqslant n\}$ in $\downarrow\!\overline{R}(F)$. It is easy to check that $G$ satisfies
(sL-1) and (sL-2) and $G\in \mathcal{S}(M, F)\ne\emptyset$. Thus $(U, R, \mathcal{F})$ is  an sL-approximation space.
\hfill$\Box$

\begin{tm} {\em (Representation Theorem for algebraic (s)L-domains)}\label{tm-asl-topsl} A poset $(L, \leqslant)$ is an algebraic sL-domain (resp., algebraic L-domain) if and only if there exists some topological sL-approximation space  (resp., topological  L-approximation space) such that $(\mathfrak{C}(U, R, \mathcal{F}), \subseteq)\cong (L, \leqslant)$.
\end{tm}
\noindent\noindent{\bf Proof.} Notice that the upper approximation operator in a topological GA-space is a closure operator. So the proof here is similar to the that of \cite[Theorems 4.2, 4.3]{Wu-Xu-FIE}. \hfill$\Box$\\
The above representation theorems demand much  in some sense. For example, to represent an sL-domain, by Theorem \ref{tm-rep-sL}, one needs use a CF-approximation space which itself is  an sL-approximation space. However, we will see that it would be sufficient if the CF-approximation space has a dense  sL-approximation subspace in the sense of Definition \ref{dn-dence-sub}.

\begin{tm}{\em(Representation Theorem by using dense subspaces)} Let $(L, \leqslant)$ be a poset. Then the following statements are equivalent.

$(1)$ $L$ is an sL-domain (resp.,  L-domain, bc-domains).

 $(2)$ There exists a dense CF-approximation subspace $(V, Q, \mathcal{G})$ of an sL-approximation space $(U, R,\mathcal{F})$ (resp.,  L-approximation space, bc-approximation subspace)  such that $(\mathfrak{C}(V, Q, \mathcal{G}), \subseteq)\cong (L, \leqslant)$.

$(3)$  There exists a CF-approximation space $(U, R,\mathcal{F})$ which has a dense  sL-approximation subspace $(V, Q, \mathcal{G})$ (resp.,  L-approximation subspace,  bc-approximation subspace) such that $(L, \leqslant)\cong(\mathfrak{C}(U, R, \mathcal{F}), \subseteq)$.
\end{tm}
\noindent{\bf Proof.}
 It directly follows from Theorems \ref{tm-rep-sL}, \ref{tm-rep-Ldo}, \ref{th-rep-bc-domain}, Proposition \ref{pn-dens-clo} and that every  CF-approximation space is a dense  CF-approximation subspace of itself.
 \hfill$\Box$\\


%

In what follows, we  first recall the  definition  of CF-approximable relations in \cite{Wu-Xu}, and then establish categorical equivalences of related categories.

\begin{dn}{(\cite{Wu-Xu})}\label{dn-CF-app}  Let $(U_{1}, R_{1}, \mathcal{F}_{1})$, $(U_{2}, R_{2}, \mathcal{F}_{2})$ be CF-approximation spaces, and $\mathrel{\Theta}\subseteq\mathcal{F}_{1}\times\mathcal{F}_{2}$ a binary relation.  If

$(1)$  for all $F\in \mathcal{F}_1$, there exists   $G\in \mathcal{F}_{2}$ such that $F\mathrel{\Theta} G$;

$(2)$  for all $F, F^{\prime}\in \mathcal{F}_{1}$, $G\in  \mathcal{F}_{2}$, if $F\subseteq \overline{R_{1}}(F^{\prime})$, $F\mathrel{\Theta} G$, then $F^{\prime}\mathrel{\Theta} G$;

$(3)$  for all $F\in \mathcal{F}_{1}$, $G, G^{\prime}\in  \mathcal{F}_{2}$, if $F\mathrel{\Theta} G$, $G^{\prime}\subseteq \overline{R_{2}}(G)$,  then $F\mathrel{\Theta} G^{\prime}$;

$(4)$ for all $F\in \mathcal{F}_{1}$,  $G\in \mathcal{F}_{2}$, if $F\mathrel{\Theta} G$, then there are $F^{\prime}\in \mathcal{F}_{1}$,  $G^{\prime}\in \mathcal{F}_{2}$ such that  $F^{\prime}\subseteq \overline{R_{1}}(F)$, $G\subseteq \overline{R_{2}}(G^{\prime})$ and $F^{\prime}\mathrel{\Theta} G^{\prime}$; and

$(5)$  for all $F\in \mathcal{F}_{1}$, $G_{1}, G_{2}\in \mathcal{F}_{2}$, if $F\mathrel{\Theta} G_{1}$ and $F\mathrel{\Theta} G_{2}$, then there exists    $G_{3}\in \mathcal{F}_{2}$ such that  $G_{1}\cup G_{2}\subseteq\overline{R_{2}}(G_{3})$ and $F\mathrel{\Theta} G_{3}$,

\noindent then $\mathrel{\Theta}$ is called a {\em CF-approximable relation} from $(U_{1}, R_{1}, \mathcal{F}_{1})$ to $(U_{2}, R_{2}, \mathcal{F}_{2})$.
\end{dn}

We denote  the catgeory of sL-approximation spaces (resp., L-approximation spaces, bc-approximation spaces, topological sL-approximation spaces, topological L-approximation spaces) by {\bf sL-APPS} (resp.,  {\bf L-APPS},  {\bf bc-APPS} {\bf TopsL-APPS}, {\bf TopL-APPS}),
  which is  a full subcategory of the categroy of CF-approximation spaces with  CF-approximable relations as morphisms.
We use  {\bf sL-DOM} (resp.,  {\bf L-DOM}, {\bf bc-DOM} {\bf  AsL-DOM}, {\bf AL-DOM}) to  denote  the category of sL-domains (resp., L-domains, bc-domains algebraic sL-domains, algebraic L-domains) with  Scott continuous mappings. By \cite[Theorem 5.9]{Wu-Xu}, we know that the category of CF-approximation spaces  is equivalent to that of continuous domains. Thus we have the following results.

\begin{tm}\label{Cat-sL-sLD}
The category {\bf sL-APPS} (resp.,  {\bf L-APPS},  {\bf bc-APPS} {\bf TopsL-APPS}, {\bf TopL-APPS}) is equivalent to {\bf sL-DOM} (resp.,  {\bf L-DOM}, {\bf bc-DOM} {\bf  AsL-DOM}, {\bf AL-DOM}).
\end{tm}

\section{Concluding remarks}

It is a common sense that representing mathematical structures by various families of subsets of other structures has proved  effective in applications and theoretic studies.  By means of  upper  approximation operators in rough set theory,   we  introduced sL-approximation spaces,  L-approximation spaces and bc-approximation spaces and gave several representation theorems by the families of CF-closed sets in such kind of spaces for sL-domains, L-domains and bc-domains. It has been proved that the collection of CF-closed sets in an sL-approximation
space (resp., an L-approximation
space, a bc-approximation space) equipped with the set-theoretic inclusion order is an sL-domain (resp., an L-domain, a bc-domain);  every sL-domain (resp., L-domain, bc-domain) is order-isomorphic to the  collection of CF-closed sets of some sL-approximation
space (resp., L-approximation spaces, bc-approximation spaces)  with respect to  the set-theoretic inclusion order. Since CF-approxiamtion spaces are generalizations of abstract bases and sL-domains are generalizations of L-domains, compare to representations of L-domains via abstract bases, our work  have a larger scope.  The work of this paper  has   both the backgrounds of domain theory and that of  rough set theory.   The  categorical
equivalence established  in this paper can be  used to further synthesize these two independent branches of mathematics  in theoretical computer science.

%




\end{document}